\documentclass[12pt]{article}
\newcommand{\eqdef}{\, =\kern -12.7pt\raise 6pt\hbox{{\tiny\textrm{def}}}\,\,}
\newcommand{\hstir}[3]{{#1\brace #2}_{\kern-2pt#3}}
\begin{document}
\newcommand{\pstir}[3]{\ensuremath{{#1\brace
#2}^{\kern-1pt\lower2pt\hbox{\scriptsize #3}}}}
\newcommand{\vanish}[1]{}
\newcommand{\obar}[1]{\overline{\omega_{#1}}}
\newtheorem{theorem}{Theorem}
\title{Some New Aspects of the Coupon-Collector's Problem}
\author{Amy N. Myers and Herbert S. Wilf\\
University of Pennsylvania\\
Philadelphia, PA 19104-6395}
\maketitle
\begin{abstract}
We extend the classical coupon collector's problem to one in which two
collectors are simultaneously and independently seeking collections of $d$ coupons. We find, 
in finite terms,
the probability that the two collectors finish at the same trial, and we
find, using the methods of Gessel-Viennot, the probability that the game
has the following ``ballot-like'' character: the two collectors are tied 
with each other
for some initial number of steps, and after that the player who first gains
the lead remains ahead throughout the game. As a by-product we obtain the evaluation in finite terms of certain infinite series whose
coefficients are powers and products of Stirling numbers of the second kind.

We study the variant of the original coupon collector's problem in which a 
single collector wants to obtain at least $h$ copies of each coupon. Here 
we give a simpler derivation of results of Newman and Shepp, and extend 
those results. Finally we obtain the distribution of the number of coupons 
that have been obtained exactly once (``singletons'') at the conclusion of 
a successful coupon collecting sequence.
\end{abstract}
\newpage
{\small\baselineskip 8pt\tableofcontents}

\bigskip

\section{Introduction and results}
The classical coupon collector's problem is the following. Suppose that a
breakfast cereal manufacturer offers a souvenir (``coupon'') hidden in each
package of cereal, and there are $d$ different kinds of souvenirs
altogether. The collector wants to have a complete collection of all $d$
souvenirs. What is the probability $p(n,d)$ that exactly $n$ boxes of
cereal will have to be purchased in order to obtain, for the first time, a
complete collection of at least one of each of the $d$ kinds of souvenir
coupons?

The answer to that question is well known (e.g., \cite{wi}, p. 132) to be
\begin{equation}
\label{eq:prob1}
p(n,d)=\frac{d!}{d^{n}}{n-1\brace d-1},
\end{equation}
where the ${n\brace k}$'s are the Stirling numbers of the second kind.

We study, in this paper, a number of other aspects of this problem, as well 
as a generalization of it to a two-player game.

First, suppose we have two coupon collectors, drawing coupons 
simultaneously, and each seeking to obtain a complete collection of $d$ 
coupons. We ask for the probability that the two games are completed at the 
same time. The answer is given by (\ref{eq:sametime}) below. That answer is 
expressed in finite terms, owing to the closed form evaluation of the 
ordinary power series generating function for the squares of the Stirling 
numbers of the second kind, contained in (\ref{eq:stirsq}).

Next we consider the following two-person game. Again two coupon collectors 
are simultaneously drawing coupons at random. This time we are interested 
in a ballot-like problem: what is the probability that the player who first 
completed a collection (the \textit{winner}) was never behind (i.e., never 
had fewer distinct coupons) at any intermediate stage of the play? Here we 
give a complete answer to a slightly easier question, namely the following: 
what is the probability that after an initial segment of play in which the 
players are tied, one of them takes the lead and keeps the lead until the 
end. The answer is in eq. (\ref{eq:fnlprob}) below, and is obtained by the 
Gessel-Viennot theory of nonintersecting lattice paths.

In each of these cases the answer can first be written as an infinite
series whose coefficients involve various products of Stirling numbers.
What is interesting, though, is that in all such cases we are able to
express the answers in finite terms. Indeed, one of our main results here is
the observation that infinite series whose coefficients involve various
powers and products of Stirling numbers of the second kind can readily be
evaluated in finite terms.

In section \ref{sec:dixie} we return to the original collecting problem of 
obtaining at least one copy of each coupon, but now we study the variant of 
the problem in which a single collector wants to obtain at least $h\ge 1$ 
copies of each coupon. We obtain the generating function (\ref{eq:hcopies}) 
for the probability that exactly $n$ trials are needed, the exact value of 
the average number of trials (\ref{eq:avgxact}), and the asymptotic 
behavior (\ref{eq:asym}) of these quantities as $n\to\infty$.

Finally, in section \ref{sec:singletons}  we study the number of coupons 
that have been collected only once, at the end of a collection sequence. We 
find the distribution function (\ref{eq:jprob}) for this number, and show 
that the average number of these singletons is just the harmonic number 
$H_d=1+1/2+\dots+1/d$.

\section{The two-person collecting competition}
\subsection{Simultaneous completion}
We find now the probability of simultaneous completion of two 
independent  coupon
collecting sequences. Evidently this is,
\begin{equation}
\sum_{n\ge 0}p(n,d)^2=\sum_{n\ge 0}\frac{d!^2}{d^{2n}}\pstir{n-1}{d-1}{2},\\
\label{eq:sametm}
\end{equation}
which expresses the answer as an infinite sum. We can rewrite this as a
finite sum by finding a finite expression for the generating function for
the squares of the Stirling numbers of the second kind,
\[F_k(x)\eqdef \sum_{n\ge k}\pstir{n}{k}{2}x^n,\]
analogously to the well known generating function for these numbers themselves,
\begin{equation}
\label{eq:pgfh1}
\sum_{n\ge k}{n\brace k}x^n=\frac{x^k}{(1-x)(1-2x)\dots (1-kx)}.
\end{equation}
The easiest way to do this is via the standard explicit formula for these 
Stirling numbers, viz.
\begin{eqnarray}
{n\brace k}&=&\frac{1}{k!}\sum_{r=1}^k(-1)^{k-r}{k\choose r}r^n\qquad (1\le
k\le n)\label{eq:stirexpl}\\
&\eqdef& \sum_{r=1}^kA_{k,r}r^{n-k},\nonumber
\end{eqnarray}
where we have written
\begin{equation}
\label{eq:aij}
A_{k,r}=\frac{(-1)^{k-r}r^k}{k!}{k\choose r}.
\end{equation}
It follows that
\begin{eqnarray}
F_k(x)\eqdef\sum_{n\ge k}\pstir{n}{k}{2}x^n&=&\sum_{n\ge
k}x^n\sum_{r,s=1}^kA_{k,r}A_{k,s}r^{n-k}s^{n-k}\nonumber\\
&=&x^k\sum_{r,s=1}^kA_{k,r}A_{k,s}\sum_{n\ge k}(rsx)^{n-k}\nonumber\\
&=&x^k\sum_{r,s=1}^k\frac{A_{k,r}A_{k,s}}{1-rsx},\qquad
(|x|<\frac{1}{k^2}).\label{eq:stirsq}
\end{eqnarray}

Thus for the simultaneous completion probability we obtain, from
(\ref{eq:sametm}),
\begin{equation}
\label{eq:sametime}
\sum_{n\ge
0}p(n,d)^2=\frac{d!^2}{d^{2d}}\sum_{r,s=1}^{d-1}\frac{A_{d-1,r}A_{d-1,s}}{1-\frac{rs}{d^
2}},
\end{equation}
by (\ref{eq:stirsq}), where the $A$'s given by (\ref{eq:aij}). This
sequence of probabilities, for $d=1,2,\dots$, begins as

\[ 1,\, \frac{1}{3},\, \frac{11}{70},\, \frac{9}{91},\,
     \frac{688877}{9561123},\, \frac{358555}{6330324},\,
\frac{2730269557627901}{58560931675094420},\,
     \frac{146271649897951}{3695016639410525},\dots ,\]
i.e., as
\[
1,0.33333..,0.15714..,0.098901..,0.072049..,0.056640..,0.046622..,0.039586..,\dots 
.\]

\subsection{Neck-and-neck, then always ahead}
We encode a sequence of $n$ draws as a path $\omega$ with $n$ vertices in the
lattice $\cal{L}$ consisting of vertices $(i,j)$, and edges
$\{(i,j),(i+1,j)\}, \{(i,j),(i+1,j+1)\}$, for all $i,j \geq 0$.  The first
coordinate of a vertex in the path gives the number of draws, or
\textit{steps}, and the second coordinate gives the number of distinct
coupons the collector has at that step.  Thus $\omega$ starts at $(0,0)$
indicating the collector has 0 coupons at draw 0, proceeds to $(1,1)$ (the
collector has 1 coupon after 1 draw), and ends at $(n,d)$, $n \geq d$ (the
collector has a complete collection at step $n$).  We write
$\omega=(0,0)\overline{\omega}(n,d)$, where $\overline{\omega}$ is a path from
$(1,1)$ to $(n-1,d-1)$, to indicate that $\omega$ starts at the vertex
$(0,0)$,
continues with the first vertex $(1,1)$ in $\overline{\omega}$, then follows
$\overline{\omega}$ through to $(n-1,d-1)$, and finally ends with the vertex
$(n,d)$.

We assign a weight of $i/d$ to each horizontal edge $\{(i,j),(i+1,j)\}$ in 
the lattice $\cal{L}$.  This is
the probability that at the $(j+1)^{\mathrm{st}}$ step, the collector draws 
one of the
$i$ distinct coupons already collected at step $j$.  We assign a weight of
$1-i/d$ to each \textit{northeast} edge $\{(i,j),(i+1,j+1)\}$.  The 
probability that the
collector draws the particular sequence of coupons encoded by the path
$\omega$ is given by the product of the weights on the edges of
$\omega$.  We let $P(\omega)$ denote this probability.

Suppose one collector, the \textit{winner}, collects all $d$ distinct coupons
for the
first time at step $n.$  (At step $n-1$ the winner had $d-1$ distinct
coupons.)
Let $\omega_1$ be the lattice path which encodes the winner's sequence of
draws.
   Let $\omega_2$ encode the other collector's draws.  We
compute the
probability $p(d)$ that $\omega_1$ and $\omega_2$ are identical until some
point at which the winner takes the lead and the other collector never
catches up.

To do this, we begin by supposing $\omega_1$ is identical to $\omega_2$ until
step $k$, at which point both collectors have $d_1$ distinct coupons. The
argument splits into two cases, namely $k\le n-2$ and $k=n-1$. In both
cases, at step
$k+1$ the winner collects one additional distinct coupon while the other
collector does not.  After step $k$, the two paths never intersect again.  The
winner collects all $d$ distinct coupons for the first time at step $n$.
Suppose the other collector has $d_2$ distinct coupons at this point.  The
probability we seek is
\begin{equation}
\label{eq:pnd}
p(d)=\sum_{n=d}^{\infty}\sum_{k=1}^{n-1}\sum_{d_1=1}^{d-1}\sum_{d_2=d_1}^{d-1}
\sum_{(\omega_1,\omega_2)} P(\omega_1)P(\omega_2)
\end{equation}
where the innermost sum ranges over all pairs $(\omega_1,\omega_2)$ described
above.

\subsection{The case $k\le n-2$}

Write $\omega_1=\alpha\obar1(n,d)$, where $\alpha$
denotes a lattice path from $(0,0)$ to $(k,d_1)$, and $\obar1$ denotes a
path from $(k+1,d_1+1)$ to $(n-1,d-1)$. Similarly, set
$\omega_2=\alpha\obar2$, where $\alpha$ is as above, and $\obar2$ is a path 
from $(k+1,d_1)$
to $(n,d_2)$.  Note that $\obar1$ and
$\obar2$ are nonintersecting paths in the lattice $\cal{L}$.
In terms of these we have $P(\omega_1)=P(\alpha)(1-d_1/d)P(\obar1)(1/d)$ and
$P(\omega_2)=P(\alpha)(d_1/d)P(\obar2)$. Hence from (\ref{eq:pnd}) we find
for the combined probability of all pairs if $k\le n-2$,
\begin{eqnarray}
\label{eq:pnd2}
p(d)_{k \leq 
n-2}&=&\sum_{n=d}^{\infty}\sum_{k=1}^{n-2}\sum_{d_1=1}^{d-1}\sum_{d_2=d_1}^{d-1}
\sum_{(\omega_1,\omega_2)}
\left[P(\alpha)\left(1-\frac{d_1}{d}\right)P(\obar1)\left(\frac{1}{d}\right)\right]\left[P(\alpha)\left(\frac{d_1}{d}\right)P(\obar2)\right]\nonumber\\
&=&\sum_{n=d}^{\infty}\sum_{k=1}^{n-2}\sum_{d_1=1}^{d-1}\sum_{d_2=d_1}^{d-1}\left(1-\frac{d_1}{d}\right)\left(\frac{1}{d}\right)\left(\frac{d_1}{d}\right)\sum_{\alpha}P(\alpha)^2\sum_{(\obar1,\obar2)}P(\obar1)P(\obar2)
\end{eqnarray}

At this point we have translated a question about coupon collecting into a 
problem
involving nonintersecting paths in a lattice.  We have set the stage for
application of the Gessel-Viennot theorem \cite{gv}.  This result concerns
pairs of
nonintersecting lattice paths
with no constraints on vertices or edges in the paths.  For this reason we
have written $\omega_1$ and $\omega_2$ in terms of $\obar1$ and
$\obar2$.

The theorem refers to an
arbitrary set $\cal{L}$, which we will take to be the lattice defined
earlier, and a weight (or valuation) $v$, which we take to be $P$.  The
theorem equates a sum of weights of paths with the determinant of a matrix
$(a_{ij})_{1 \leq i,j \leq l}$.  The entries of this matrix are defined by
$a_{ij}=\sum_{\omega}v(\omega)$, where $\omega$ ranges over all paths from
$A_i$ to $B_j$.

The theorem requires that two given sequences, $(A_1,A_2,\dots,A_l)$ and
$(B_1,B_2,\dots,B_l)$, of vertices in $\cal{L}$, the sets $\Omega_{ij}$,
$1\leq i,j \leq l$, of all paths in $\cal{L}$ between $A_i$ and $B_j$, and
the weight $v$ satisfy both the \textit{finiteness} and \textit{crossing
conditions}.  The \textit{finiteness condition} requires the set of paths
in $\Omega_{ij}$ with nonzero weight be finite.  The \textit{crossing
condition} requires that paths in $\Omega_{ij'}$ and $\Omega_{i'j}$, $i<i'$
and $j<j'$, with nonzero weight share a common vertex.  Both conditions
hold for the paths we consider.\\

\begin{theorem} (Gessel-Viennot)  Suppose $\cal{L}$, $v$,
$(A_1,A_2,\dots,A_l)$, and $(B_1,B_2,\dots,B_l)$ satisfy both the finiteness
and crossing conditions.  Then the determinant of the matrix $(a_{ij})_{1
\leq i,j \leq l}$ is the sum of the weights of all configurations of paths
$(\omega_1,\omega_2,\dots,\omega_l)$ satisfying the following two conditions:\\
\indent (i)  The paths $\omega_k$ are pairwise nonintersecting, and\\
\indent (ii)  $\omega_k$ is a path from $A_k$ to $B_k$.\\
In other words,
\[
\det\left(\left\{a_{ij}\right\}_{i,j=1}^l\right) =
\sum_{(\omega_1,\omega_2,\dots,\omega_l)} v(\omega_1)v(\omega_2) \dots
v(\omega_l).
\]
\end{theorem}

Application of this theorem to our problem requires the computation of only
a $2
\times 2$ determinant!  Let $A_1=(k+1,d_1+1)$, $A_2=(k+1,d_1)$,
$B_1=(n-1,d-1)$,
and $B_2=(n,d_2)$.  Then
\begin{equation}
\label{eq:tot}
\sum_{(\obar1,\obar2)}P(\obar1)P(\obar2) = \det
\left[
\begin{array}{ll}
a_{11}&a_{12}\\
a_{21}&a_{22}
\end{array}
\right]
\end{equation}
where $a_{ij}$ is the sum $\sum_{\omega}P(\omega)$ over all paths $\omega$ 
from $A_i$ to $B_j$.

\subsection{Paths from $A$ to $B$}
\label{sec:AtoB}

In this section, we compute the probability $P(\omega)$ of an arbitrary 
path $\omega$ from a vertex $A=(a_1,b_1)$ to a vertex $B=(a_2,b_2)$, as well as the sum over all such paths.  Such 
a path contains $b_2-b_1$ \textit{northeast} edges $\{(i,j),(i+1,j+1)\}$ 
and $(a_2-a_1)-(b_2-b_1)$  horizontal edges.  The weights
assigned to northeast edges in order from left to right are
$1-\frac{b_1}{d},1-\frac{b_1+1}{d},\dots,1-\frac{b_2-1}{d}$.  The weight 
assigned to a horizontal edge depends its coordinates.  Consider the edge
$\{(i,j),(i+1,j)\}$.  This edge indicates the collector has $j$ distinct
coupons at step $i$ and draws one of the same $j$ coupons at step
$i+1$.  The probability of this (weight of the edge) is
$\frac{j}{d}$.  Thus the probability of a path $\omega$ from $A$ to $B$ is
\begin{eqnarray}
\label{eq:AtoB}
P(\omega)&=&
\left(\frac{b_1}{d}\right)^{e_1}\left(1-\frac{b_1}{d}\right)\left(\frac{b_1+1}{d}\right)^{e_2}\left(1-\frac{b_1+1}{d}\right)\dots\left(1-\frac{b_2-1}{d}\right)\left(\frac{b_2}{d}\right)^{e_{b_2-b_1+1}
}\nonumber\\
&=&\frac{1}{d^{a_2-a_1}}\frac{(d-b_1)!}{(d-b_2)!}
(b_1)^{e_1}(b_1+1)^{e_2}\dots(b_2)^{e_{b_2-b_1+1}}
\end{eqnarray}
where $e=(e_1,e_2,\dots,e_{b_2-b_1+1})$ is an ordered partition, a
\textit{composition}, of $(a_2-a_1)-(b_2-b_1)$ into $b_2-b_1+1$ nonnegative 
integer parts.  With this we compute the sum of the probabilities of all 
paths from $A$ to $B$.
\begin{eqnarray}
\label{eq:AtoBsum}
\sum_{\omega=A\cdots B}P(\omega)&=&\sum_e
\left(\frac{b_1}{d}\right)^{e_1}\left(1-\frac{b_1}{d}\right)\left(\frac{b_1+1}{d}\right)^{e_2}\left(1-\frac{b_1+1}{d}\right)\dots\left(1-\frac{b_2-1}{d}\right)\left(\frac{b_2}{d}\right)^{e_{b_2-b_1+1}
}\nonumber\\
&=&\frac{1}{d^{a_2-a_1}}\frac{(d-b_1)!}{(d-b_2)!}\sum_e
(b_1)^{e_1}(b_1+1)^{e_2}\dots(b_2)^{e_{b_2-b_1+1}}
\end{eqnarray}
where the sum is over all compositions $e=(e_1,e_2,\dots,d_{b_2-b_1+1})$ of 
$(a_2-a_1)-(b_2-b_1)$ into $b_2-b_1+1$ nonnegative integer parts.
This is the
coefficient of $x^{(a_2-a_1)-(b_2-b_1)}$ in the series expansion of
\[\frac{1}{(1-b_1x)(1-(b_1+1)x)\dots (1-b_2x)}\]
so we can find a simpler formula for it by looking at the partial fraction
expansion
\begin{equation}
\label{eq:parfrac}
\frac{1}{\prod_{m=a}^{b}(1-mx)}=\sum_{m=a}^{b}\frac{B_m}{1-mx},
\end{equation}
where
$$
B_m=\frac{(-1)^{b-m}m^{b-a}}{(b-a)!}{b-a\choose m-a}.
$$
 From this and (\ref{eq:AtoBsum}) we obtain
\begin{eqnarray}
\sum_{\omega=A\cdots B} P(\omega)&=&
\frac{1}{d^{a_2-a_1}}\frac{(d-b_1)!}{(d-b_2)!}[x^{(a_2-a_1)-(b_2-b_1)}]\left\{\sum_{m=b_1}^{b_2}\frac{B_m}{1-mx}\right\}\nonumber\\
&=&\frac{1}{d^{a_2-a_1}}\frac{(d-b_1)!}{(d-b_2)!}\sum_{m=b_1}^{b_2}B_mm^{(a_2-a_1)-(b_2-b_1)}\nonumber\\
&=&\frac{1}{d^{a_2-a_1}}\,\frac{(d-b_1)!}{(d-b_2)!}\,\sum_{m=b_1}^{b_2}\frac{(-1)^{b_2-m}m^{a_2-a_1}}{(b_2-b_1)!}{b_2-b_1 
\choose m-b_1}\label{eq:AtoBf}
\end{eqnarray}

\subsection{Evaluating the determinant}

We use the results of the previous section to evaluate the determinant in 
(\ref{eq:tot}).  To compute $a_{11}$, we substitute $A=A_1=(k+1,d_1+1)$ and 
$B=B_1=(n-1,d-1)$ in (\ref{eq:AtoBf}).  This yields
\begin{eqnarray}
a_{11}&=&\frac{1}{d^{n-k-2}}\,(d-d_1-1)!\,\sum_{m=d_1+1}^{d-1}\frac{{{\left( 
-1 \right)}^{d-m-1}}\,m^{n-k-2}}{(d-d_1-2)!}{d-d_1-2 \choose m-d_1-1} 
\label{eq:a11f}
\end{eqnarray}
In a similar manner we obtain
\begin{eqnarray}
a_{12}&=&\frac{1}{d^{n-k-1}}\,\frac{(d-d_1-1)!}{(d-d_2)!}\,\sum_{m=d_1+1}^{d_2}\frac{{{\left( 
-1 \right)}^{d_2-m}}\,m^{n-k-1}}{(d_2-d_1-1)!} {d_2-d_1-1 \choose 
m-d_1-1}\label{eq:a12f}\\
a_{21}&=&\frac{1}{d^{n-k-2}}\,(d-d_1)!\,\sum_{m=d_1}^{d-1}\frac{{{\left( -1 
\right)}^{d-m-1}}\,m^{n-k-2}}{(d-d_1-1)!} {d-d_1-1 \choose 
m-d_1}\label{eq:a21f}\\
a_{22}&=&\frac{1}{d^{n-k-1}}\,\frac{(d-d_1)!}{(d-d_2)!}\,\sum_{m=d_1}^{d_2}\frac{{{\left( 
-1 \right)}^{d_2-m}}\,m^{n-k-1}}{(d_2-d_1)!}{d_2-d_1 \choose 
m-d_1}\label{eq:a22f}
\end{eqnarray}

Using (\ref{eq:a11f})-(\ref{eq:a22f}) we compute the determinant of
our $2\times 2$ matrix.
\begin{eqnarray}
\det
\left[
\begin{array}{ll}
a_{11}&a_{12}\\
a_{21}&a_{22}\\
\end{array}
\right]&=&
\frac{(d-d_1)!(d-d_1-1)!}{d^{2n-2k-3}(d-d_2)!}\times \\
&&\sum_{l=d_1}^{d}\sum_{m=d_1}^{d_2}\frac{m(lm)^{n-k-2}(-1)^{d_1+d_2}(l-d)(m^2-l^2)}{(d-d_1-1)!(d_2-d_1)!}{d-d_1-1 
\choose l-d_1}{d_2-d_1 \choose m-d_1}\nonumber\\
&\eqdef&\det(d,d_1,d_2,k,n)\label{eq:alphdet}
\end{eqnarray}
Substituting (\ref{eq:alphdet}) in (\ref{eq:tot}), we obtain
\begin{equation}
\label{eq:opaths}
\sum_{(\obar1,\obar2)}P(\obar1)P(\obar2)=\det(d,d_1,d_2,k,n)
\end{equation}

\subsection{The initial common segment}
\label{sec:initseg}

In the previous section we evaluated the determinant in (\ref{eq:tot}).  In 
this section we compute the sum $\sum_{\alpha} P(\alpha)^2$ in 
(\ref{eq:pnd2}).  Recall $\alpha$ is a path from $(0,0)$ to $(k,d_1)$.

Equation (\ref{eq:AtoB}) gives the probability of an arbitrary path from 
$A$ to $B$.  Substituting $A=(0,0)$ and $B=(k,d_1)$ gives the probability
$$P(\alpha)=\frac{d!}{d^k(d-d_1)!}1^{e_1}2^{e_2}\cdots {d_1}^{e_{d_1}}$$
of an arbitrary path $\alpha$ from $(0,0)$ to $(k,d_1)$.  It follows that
\begin{eqnarray}
\label{eq:tailAB}
\sum_{\alpha=(0,0)\cdots(k,d_1)}P(\alpha)^2&=&\frac{d!^2}{d^{2k}(d-d_1)!^2}\sum_{e_1+\dots+e_{d_1}=k-d_1}1^{2e_1}2^{2e_2}\dots d_1^{2e_{d_1}}\nonumber\\
&=&\frac{d!^2}{d^{2k}(d-d_1)!^2}[x^{k-d_1}]\left\{\frac{1}{(1-1^2x)(1-2^2x)\dots (1-d_1^2x)}\right\}\nonumber\\
&=&\frac{d!^2}{d^{2k}(d-d_1)!^2}\sum_{m=1}^{d_1}C_mm^{2k-2d_1}\qquad 
(\mathrm{as\ in}\ (\ref{eq:parfrac}))\nonumber\\
&=&\frac{2d!^2}{(d-d_1)!^2d^{2k}(2d_1)!}\sum_{m\ge 
1}(-1)^{d_1-m}{2d_1\choose d_1+m}m^{2k}.\label{eq:sumsq}\\
&\eqdef&\mathrm{init}(d,d_1,k)\label{eq:initseg}
\end{eqnarray}

\subsection{The case $k=n-1$}

Suppose now that the two walks are identical up to the point $(n-1,d-1)$.
Since step $n$ is the finish, the next step for the winning player will be
to $(n,d)$, and for the loser, to $(n,d-1)$. These last steps have
respective probabilities $1/d$ and $1-1/d$. Hence the probability of the
complete pair of walks in this case is the probability of two identical
walks from $(0,0)$ to $(n-1,d-1)$ (which is given by (\ref{eq:sumsq}) with
$(k,d_1):=(n-1,d-1)$) multiplied by $(d-1)/d^2$.

\subsection{Putting it together}

We now substitute (\ref{eq:alphdet}) and (\ref{eq:initseg}) into 
(\ref{eq:pnd2}) to obtain the probability of all pairs
of paths that we are considering,
\begin{eqnarray}
p(d)&=&\sum_{n=d}^{\infty}\sum_{k=1}^{n-2}\sum_{d_1=1}^{d-1}\sum_{d_2=d_1}^{d-1}\left(1-\frac{d_1}{d}\right)\left(\frac{d_1}{d}\right)\left(\frac{1}{d}\right)\mathrm{init}(d,d_1,k)\mathrm{det}(d,d_1,d_2,k,n)\nonumber\\
&&\qquad\qquad+\frac{d-1}{d^2}\sum_{n=d}^{\infty}\textrm{init}(d,d-1,n-1),
\label{eq:totalprob}\\
&\eqdef&\Sigma_1+\Sigma_2\label{Sigmas}
\end{eqnarray}

It turns out that the sums over the indices $d_2,n,k$ can all be carried
out in explicit closed form. Hence we can obtain an expression  which is in
finite terms for the total probability.

First, the sum on $d_2$ in $\Sigma_1$ above can be done in closed form since
\[\sum_{d_2=d_1}^{d-1}(-1)^{d_2}{d-d_1\choose d-d_2}{d_2-d_1\choose
t-d_1}=(-1)^{d+1}{d-d_1\choose d-t}.\]
Next, the remaining sum over the indices $n$ and $k$, in the first
summation, $\Sigma_1$ is
\begin{equation}
\label{eq:psidef}
\psi(d,r,s,t)\eqdef\sum_{n=d}^{\infty}\,\sum_{k=1}^{n-2}\frac{r^{2k}t^{n-k-1}s^{n-k-2}}{d^{2n}}=\cases{\frac{r^{2d}(d^3-2d^2-r^2d+3r^2)}{d^{2d-2}(d^2-r^2)^2s^2t},&if 
$r^2=st$;\cr
\frac{r^4(st)^{d-1}(r^2-d^2)+str^{2d}(d^2-st)}{d^{2d-2}(d^2-st)(d^2-r^2)(r^2-st)r^2s},&otherwise.}
\end{equation}
The sum over $n$ in $\Sigma_2$ is trivial, and so there remain no infinite
sums in our final expression for the probability $p(d)$, which is
\begin{eqnarray}
\label{eq:fnlprob}
&&\kern-10pt\sum_{d_1=1}^{d-1}\frac{2d!^2d_1(d-d_1)}{(d-d_1)!^2(2d_1)!}\sum_{r,s,t\ge 
1}(-1)^{d_1-r-s-t}(s-t){2d_1\choose d_1+r}{d-d_1-1\choose
s-d_1}{d-d_1\choose d-t}\psi(d,r,s,t)\nonumber\\
&&\qquad\qquad
+\frac{4(d-1)d!^2}{d^{2d-2}(2d-2)!}\sum_{r=1}^{d-1}(-1)^{d-1-r}{2d-2\choose
d-1+r}\frac{r^{2d-2}}{d^2-r^2}+\delta_{d,1}
\end{eqnarray}
where $\psi$ is given by (\ref{eq:psidef}).

This is the probability that the game is of the type we described, namely
where the players are tied for some initial segment of trials and then the
player who pulls ahead remains ahead always, expressed as a finite sum
(albeit a complicated one!). More precisely, the values of $p(d)$ can be
calculated, as rational numbers, with $O(d^4)$ evaluations of the above
summand. The exact values of $p(d)$, for $d=1,2,3,4,5,\dots$ are
\[\left\{
1,{\frac{2}{3}},{\frac{43}{70}},{\frac{986}{2275}},{\frac{5672893}{1912246}},\dots.\right\}\]
As decimals, the values of $\{p(d)\}_{d=1}^{10}$ are
\[\left\{1.0,
0.66667,0.61429,0.43341,0.29667,0.21177,0.16016,0.12748,0.10551,0.08988\right\}.\]

\subsection{One collector never behind}

In contrast to the problem of staying ahead as soon as the tie is broken, which we have solved in the preceding sections, the problem in which the ultimate winner has never been behind is unsolved.

Suppose the winner collects all $d$ distinct coupons for the first time at
step $n$, at which point the other collector has $d'<d$ distinct 
coupons.  We discuss the probability $b(d)$ the winner has never been 
behind.  We use $b$ for ``ballot'' since this version of the problem has a 
distinct \textit{ballot-problem} flavor (see \cite{}).

Let $w_1$ be the lattice path which encodes the winner's sequence of 
draws.  Let $\omega_2$ encode the other collector's sequence of 
draws.  Then $b(d)$ is the probability that $\omega_2$ \textit{does not 
cross} $\omega_1$.  To say $\omega_2$ \textit{does not cross} $\omega_1$ 
means for each horizontal coordinate $i$ shared by vertices $(i,j_1)$ in 
$\omega_1$ and $(i,j_2)$ in $\omega_2$, we have $j_2 \leq j_1$.  In the 
case $j_2=j_1$, we say $\omega_1$ and $\omega_2$ \textit{intersect} at 
$(i,j_1)=(i,j_2)$.  Thus we seek all pairs $(\omega_1,\omega_2)$ such that 
$\omega_1$ is a path from $(0,0)$ to $(n,d)$ including the vertex 
$(n-1,d-1)$, $\omega_2$ is a path from $(0,0)$ to $(n,d')$ for $1 \leq d' 
\leq d$, and $\omega_2$ does not cross $\omega_1$.  Such a pair 
$(\omega_1,\omega_2)$ is illustrated by Figure \ref{fig:winner}.  Note that 
$\omega_1$ and $\omega_2$ may intersect several times.  The probability we 
seek is
$$b(d)=\sum_{n=d}^{\infty}\sum_{d'=1}^{d-1}\sum_{(\omega_1,\omega_2)} 
P(\omega_1)P(\omega_2)$$
where the innermost sum ranges over all pairs described above.

\setlength{\unitlength}{.2in}
\thicklines
\begin{center}
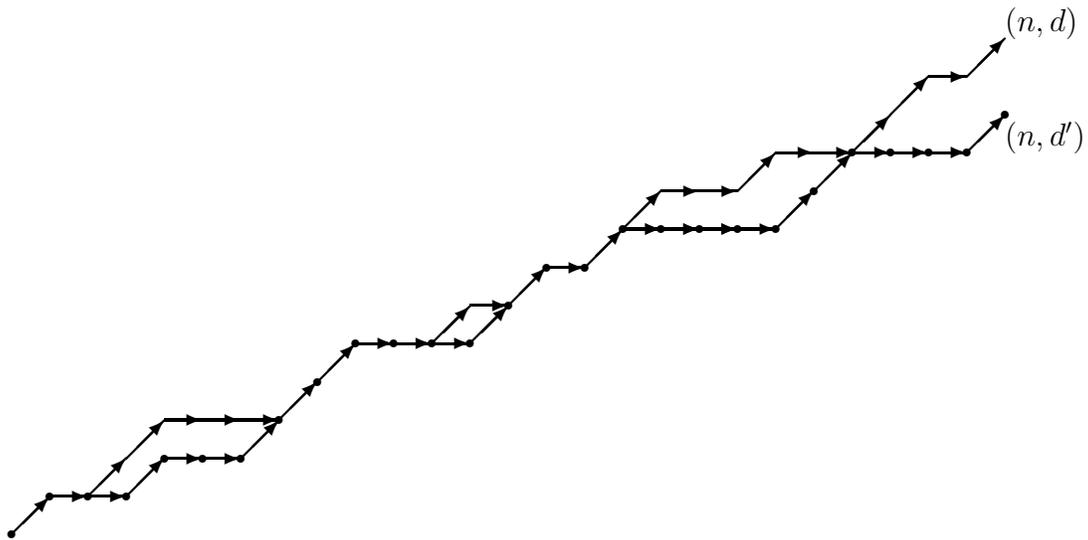
\begin{figure}
\begin{picture}(26,12)(-2,0)
\put(0,0){\circle*{.2}}
\put(1,1){\circle*{.2}}
\put(2,1){\circle*{.2}}
\put(3,1){\circle*{.2}}
\put(4,2){\circle*{.2}}
\put(5,2){\circle*{.2}}
\put(6,2){\circle*{.2}}
\put(7,3){\circle*{.2}}
\put(8,4){\circle*{.2}}
\put(9,5){\circle*{.2}}
\put(10,5){\circle*{.2}}
\put(11,5){\circle*{.2}}
\put(12,5){\circle*{.2}}
\put(13,6){\circle*{.2}}
\put(14,7){\circle*{.2}}
\put(15,7){\circle*{.2}}
\put(16,8){\circle*{.2}}
\put(17,8){\circle*{.2}}
\put(18,8){\circle*{.2}}
\put(19,8){\circle*{.2}}
\put(20,8){\circle*{.2}}
\put(21,9){\circle*{.2}}
\put(22,10){\circle*{.2}}
\put(23,10){\circle*{.2}}
\put(24,10){\circle*{.2}}
\put(25,10){\circle*{.2}}
\put(26,11){\circle*{.2}}
\put(0,0){\vector(1,1){1}}
\put(1,1){\vector(1,0){1}}
\put(2,1){\vector(1,0){1}}
\put(3,1){\vector(1,1){1}}
\put(4,2){\vector(1,0){1}}
\put(5,2){\vector(1,0){1}}
\put(6,2){\vector(1,1){1}}
\put(7,3){\vector(1,1){1}}
\put(8,4){\vector(1,1){1}}
\put(9,5){\vector(1,0){1}}
\put(10,5){\vector(1,0){1}}
\put(11,5){\vector(1,0){1}}
\put(12,5){\vector(1,1){1}}
\put(13,6){\vector(1,1){1}}
\put(14,7){\vector(1,0){1}}
\put(15,7){\vector(1,1){1}}
\put(16,8){\vector(1,0){1}}
\put(17,8){\vector(1,0){1}}
\put(18,8){\vector(1,0){1}}
\put(19,8){\vector(1,0){1}}
\put(20,8){\vector(1,1){1}}
\put(21,9){\vector(1,1){1}}
\put(22,10){\vector(1,0){1}}
\put(23,10){\vector(1,0){1}}
\put(24,10){\vector(1,0){1}}
\put(25,10){\vector(1,1){1}}
\put(2,1){\vector(1,1){1}}
\put(3,2){\vector(1,1){1}}
\put(4,3){\vector(1,0){1}}
\put(5,3){\vector(1,0){1}}
\put(6,3){\vector(1,0){1}}
\put(11,5){\vector(1,1){1}}
\put(12,6){\vector(1,0){1}}
\put(16,8){\vector(1,1){1}}
\put(17,9){\vector(1,0){1}}
\put(18,9){\vector(1,0){1}}
\put(16,8){\vector(1,0){1}}
\put(17,8){\vector(1,0){1}}
\put(18,8){\vector(1,0){1}}
\put(19,9){\vector(1,1){1}}
\put(20,10){\vector(1,0){1}}
\put(21,10){\vector(1,0){1}}
\put(22,10){\vector(1,1){1}}
\put(23,11){\vector(1,1){1}}
\put(24,12){\vector(1,0){1}}
\put(25,12){\vector(1,1){1}}
\put(26,13.2){$(n,d)$}
\put(26,10.2){$(n,d')$}
\end{picture}
\caption{The winner is never behind}
\label{fig:winner}
\end{figure}
\end{center}

Look again at Figure \ref{fig:winner}.  A pair $(\omega_1,\omega_2)$ 
appears to form a chain of flying kites anchored to the ground at 
$(0,0)$.  The highest kite has two ribbons attached to its tip.  Their 
loose ends are at $(n,d)$ and $(n,d')$.

\setlength{\unitlength}{.3in}
\thicklines
\begin{center}
\begin{figure}
\begin{picture}(13.3,12)(-4,-9.5)
\put(0,0){\circle*{.2}}
\put(1,1){\circle*{.2}}
\put(2,2){\circle*{.2}}
\put(3,2){\circle*{.2}}
\put(4,3){\circle*{.2}}
\put(5,3){\circle*{.2}}
\put(0,0){\vector(1,1){1}}
\put(1,1){\vector(1,1){1}}
\put(2,2){\vector(1,0){1}}
\put(3,2){\vector(1,1){1}}
\put(4,3){\vector(1,0){1}}
\put(-.8,-.6){\small $\mathbf{(i_0,j_0)}$}
\put(5.3,3){\small $\mathbf{(i_1,j_1)}$}
\put(8.2,-.5){\small $\mathbf{(i_1,j_1)}$}
\put(12.3,2.4){\small $\mathbf{(i_2,j_2)}$}
\put(9,0){\circle*{.2}}
\put(10,0){\circle*{.2}}
\put(11,0){\circle*{.2}}
\put(10,1){\circle*{.2}}
\put(11,2){\circle*{.2}}
\put(12,1){\circle*{.2}}
\put(12,2){\vector(1,0){1}}
\put(9,0){\vector(1,0){1}}
\put(9,0){\vector(1,1){1}}
\put(10,0){\vector(1,0){1}}
\put(11,0){\vector(1,1){1}}
\put(10,1){\vector(1,1){1}}
\put(11,2){\vector(1,0){1}}
\put(12,1){\vector(1,1){1}}
\put(12,2){\vector(1,0){1}}
\put(1.6,-2){\textit{A tail, plus ..\hspace{1.4in} .. a frame, equals ..}}
\put(1,-9){\circle*{.2}}
\put(2,-8){\circle*{.2}}
\put(3,-7){\circle*{.2}}
\put(4,-7){\circle*{.2}}
\put(5,-6){\circle*{.2}}
\put(6,-6){\circle*{.2}}
\put(1,-9){\vector(1,1){1}}
\put(2,-8){\vector(1,1){1}}
\put(3,-7){\vector(1,0){1}}
\put(4,-7){\vector(1,1){1}}
\put(5,-6){\vector(1,0){1}}
\put(5.2,-6.6){\small $\mathbf{(i_1,j_1)}$}
\put(10.2,-4){\small $\mathbf{(i_2,j_2)}$}
\put(.3,-9.5){\small $\mathbf{(i_0,j_0)}$}
\put(6,-6){\circle*{.2}}
\put(7,-6){\circle*{.2}}
\put(8,-6){\circle*{.2}}
\put(7,-5){\circle*{.2}}
\put(8,-4){\circle*{.2}}
\put(9,-5){\circle*{.2}}
\put(9,-4){\vector(1,0){1}}
\put(6,-6){\vector(1,0){1}}
\put(6,-6){\vector(1,1){1}}
\put(7,-6){\vector(1,0){1}}
\put(8,-6){\vector(1,1){1}}
\put(7,-5){\vector(1,1){1}}
\put(8,-4){\vector(1,0){1}}
\put(9,-5){\vector(1,1){1}}
\put(9,-4){\vector(1,0){1}}
\end{picture}
\caption{A kite}
\label{fig:kite}
\end{figure}
\end{center}

Each kite consists of a \textit{frame} together with a \textit{tail}.  See 
Figure \ref{fig:kite}.  A \textit{frame} from $(i_1,j_1)$ to $(i_2,j_2)$ 
consists of a pair of paths from $(i_1,j_1)$, the \textit{lower tip} of the 
frame, to $(i_2,j_2)$, the \textit{upper tip}, which intersect only at the 
endpoints.  A \textit{tail} from $(i_1,j_1)$ to $(i_2,j_2)$ consists of two 
identical paths between these endpoints.  The \textit{length} of a tail is 
the number of vertices in the tail minus one, i.e., the number of edges.

A pair $(\omega_1,\omega_2)$ such that $\omega_2$ does not cross $\omega_1$ 
forms an alternating sequence of tails and frames, beginning with a 
tail.  Note that tails may have length zero.

The upper tip of the final frame in this sequence is the common endpoint 
for two paths which intersect only at this common endpoint (these are the 
``ribbons'' described above).  One path ends at $(n,d)$, this is the 
\textit{top ribbon}, and the other ends at $(n,d')$, the \textit{bottom 
ribbon}.

Below we compute the probability of a frame from $(i_1,j_1)$ to 
$(i_2,j_2)$, a tail from $(i_1,j_1)$ to $(i_2,j_2)$, and a pair of ribbons 
with common initial point $(k,d'')$ and terminal points at $(n,d)$ and 
$(n,d')$, respectively.

Let $f_{(i_1,j_1)}^{(i_2,j_2)}(d)$ denote the probability of a frame from 
$(i_1,j_1)$ to $(i_2,j_2)$.  Note for $f_{(i_1,j_1)}^{(i_2,j_2)}(d) \neq 
0$, we must have $i_2 \geq i_1+2$, $j_2>j_1$, and $j_2-j_1 \leq 
i_2-i_1-1$.  Assuming these conditions, we write
\begin{equation}
\label{eq:frameab}
f_{(i_1,j_1)}^{(i_2,j_2)}(d)=\sum_{(\alpha,\beta)}P(\alpha)P(\beta)
\end{equation}
where $(\alpha,\beta)$ is a pair of paths from $(i_1,j_1)$ to $(i_2,j_2)$ 
intersecting only at the endpoints such that $\beta$ does not cross 
$\alpha$ (i.e., $\alpha$ forms the upper edge of the frame, and $\beta$ 
forms the lower edge).

We convert the sum above into a determinant using the Gessel-Viennot 
thereom.  Evaluation of the determinant gives
\[
f_{(i_1,j_1)}^{(i_2,j_2)}(d)=\frac{j_1 j_2 (d-j_1)!^2}{d^{2(i_2-i_1)} 
(j_2-j_1)!^2 (d-j_2)!^2}\sum_{l,m=j_1}^{j_2} 
(-1)^{l+m}(lm)^{i_2-i_1-2}(l-j_1)(m-j_2){j_2-j_1 \choose m-j_1} {j_2-j_1 
\choose l-j_1}
\]

We compute the probability $t_{(i_1,j_1)}^{(i_2,j_2)}(d)$ of a tail from 
$(i_1,j_1)$ to $(i_2,j_2)$ in a manner analogous to the computation of 
$\sum_{\alpha}P(\alpha)^2$ in section \ref{sec:initseg}.  We obtain
\[
t_{(i_1,j_1)}^{(i_2,j_2)}(d)=\frac{(d-j_1)!^2}{d^{2(i_2-i_1)}(2j_2)!(d-j_2)!^2 
} \sum_{m=1}^{j_2} (-1)^{j_2-j_1} m^{2(i_2-i_1)} (2m)! {2j_2 \choose j_2+m} 
{m+j_1-1 \choose j_1-m}
\]

Finally we compute the probability $r(d,d',d'',k,n)$ of a pair of ribbons 
with common initial point $(k,d'')$ and terminal points $(n,d')$ and 
$(n,d)$.  The probability is given by a determinant similar to the one in 
(\ref{eq:tot}).  In the present case, we have $d''$ in place of $d_1$ and 
$d'$ in place of $d_2$.  Thus
\[
r(d,d',d'',k,n)=\det(d,d',d'',k,n)
\]

\vanish{
\section{Winning margin}

Now we look for the probability distribution of the number of distinct
coupons that the second player has collected at the moment the first player
completes the collection. Let $g(d,d')$ denote the probability that the
second player has collected exactly $d'$ distinct coupons at that moment.
The probability that the first player finishes after exactly $n$ trials is
$p(n,d)$, of eq. (\ref{eq:prob1}). The probability that the second player
has exactly $d'$ distinct coupons after $n$ trials, given that the first
player has just completed the collection at that time, is
\[{d\choose d'}{n\brace d'}\frac{d'!}{d^n},\]
if $1\le d'<d$.
Thus for $d'<d$ our distribution $g$ is given by
\begin{eqnarray*}
g(d,d')&=&\sum_{n\ge d'}\frac{d!}{d^{n}}{n-1\brace d-1}{d\choose
d'}{n\brace d'}\frac{d'!}{d^n}\\
&=&\frac{d!^2}{(d-d')!}\sum_{n\ge d'}{n-1\brace d-1}{n\brace
d'}\frac{1}{d^{2n}}\\
&=&\frac{d!^2}{(d-d')!d^{2d'}}\sum_{r=1}^{d-1}r^{d'-d}\sum_{s=1}^{d'}\frac{A_{d-
1,r}A_{d',s}}{1-\frac{rs}{d^2}}-\delta_{d',1},
\end{eqnarray*}
by (\ref{eq:stirprod}). A table of the probabilities $\{g(6,j)\}_{j=1}^5$
is as follows:
\[ .000003  , .000793, .018444, .118454, .333986 \]
These do not sum to 1 because the second player might have completed a
collection at some time before the first player did.}
\section{The ``double dixie-cup problem,'' of Newman and Shepp, revisited}
\label{sec:dixie}
Here we consider a different generalization of the coupon collector's problem. Let 
integers $h, d\ge 1$ be fixed. Again we are sampling with replacement from 
$d$ kinds of coupons, but now $T$ is the epoch at which we have collected 
at least $h$ copies of each of the $d$ coupons, for the first time (for 
example, my $h-1$ siblings and I might each want to have our own copy of 
every one of the available baseball cards). We study the expectation, the 
probability generating function, and the asymptotic behavior of the 
expectation, of this generalized problem.

These questions were investigated by Newman and Shepp \cite{ns} and the 
asymptotics were refined by Erd\H os-R\' enyi \cite{er}. It is interesting to note that this problem is equivalent to one about the evolution of a random graph. Suppose we fix $n$ vertices, and then we begin to collect from among $n$ kinds of coupons. If we collect a particular sequence, say, $\{c_1,c_2,c_3,\dots\}$ then we add the edges $(c_1,c_2)$,$(c_3,c_4)\dots$. That is, we add an edge each time we choose a new pair of coupons. Our problem about collecting at least $h$ copies of each kind of coupon is thereby equivalent to the question of obtaining a minimum degree of at least $h$ in an evolving random graph.\footnote{Our thanks to Ed Bender and to a helpful referee for pointing this out.} In this section we will not add anything new to the asymptotics of this problem. Instead we claim only a simpler derivation than the original, and an explicit generating function, which gives a nice road to the asymptotics. We deal only with generating functions in one variable, whereas in \cite{ns} multivariate generating functions were used. We obtain not only the expectation of the time to reach a collection that has at least $h$ copies of each kind of coupon, but also the complete probability distribution of that time.

For $n$ fixed, consider a sequence of $n$ drawings of coupons that 
constitutes, for the first time at the $n$th drawing, a complete collection 
of at least $h$ copies of each of the $d$ kinds of coupons.

There are $d$ possibilities for the coupon that completes the collection on 
the $n$th drawing. There are ${n-1\choose h-1}$ ways to choose the set of 
earlier drawings on which that last coupon type occurred. On the remaining 
$n-h$ drawings we can define, as usual, an equivalence relation: two 
drawings $i,j$ are equivalent if the same kind of coupon was drawn at the 
$i$th and the $j$th drawings. The number of such equivalence relations is 
equal to the number of ordered partitions of a set of $n-h$ elements into 
$d-1$ classes, each class containing at least $h$ elements. We will denote 
this latter number by $(d-1)!\hstir{n-h}{d-1}{h}$, where the 
$\hstir{n}{k}{h}$'s count the unordered partitions of an $n$-set into $k$ 
classes of at least $h$ elements each.

The number of sequences of $n$ drawings for which we achieve a complete 
collection for the first time at the $n$th drawing is therefore
\[d{n-1\choose h-1}(d-1)!\hstir{n-h}{d-1}{h}.\]
Since there are $d^n$ possible drawing sequences of length $n$, the 
probability that $T=n$ is
\begin{equation}
p_n=\frac{d!}{d^n}{n-1\choose h-1}\hstir{n-h}{d-1}{h},
\end{equation}
and the probability generating function is
\begin{eqnarray}
P_h(x)&\eqdef&\sum_{n\ge 0}p_nx^n=\sum_{n\ge 0}\frac{d!}{d^n}{n-1\choose 
h-1}\hstir{n-h}{d-1}{h}x^n\nonumber\\
\label{eq:pgf}
&=&d!{xD-1\choose h-1}\sum_{n\ge 
0}\hstir{n-h}{d-1}{h}\left(\frac{x}{d}\right)^n,
\end{eqnarray}
where $D=\partial/\partial x$.

It remains to find the ordinary power series generating function of the 
$\hstir{n}{k}{h}$'s. The exponential formula gives us immediately their 
exponential generating function, as
\begin{equation}
\sum_{n\ge 0}\hstir{n}{k}{h}\frac{x^n}{n!}=\frac{1}{k!}\left(e^x-1-x-\dots 
-\frac{x^{h-1}}{(h-1)!}\right)^k
\end{equation}
We can convert this into an ordinary power series generating function by 
applying the Laplace transform operator
\[\int_0^{\infty}e^{-sx}\cdots dx\]
to both sides, which yields
\[\sum_{n\ge 
0}\hstir{n}{k}{h}\frac{1}{s^{n+1}}=\frac{1}{k!}\int_0^{\infty}e^{-sx}\left(e^x-1-x-\dots 
-\frac{x^{h-1}}{(h-1)!}\right)^kdx,\]
or finally
\begin{equation}
\label{eq:ptnsgf}
\sum_{n\ge 
0}\hstir{n}{k}{h}t^{n}=\frac{1}{k!t}\int_0^{\infty}e^{-x/t}\left(e^x-1-x-\dots 
-\frac{x^{h-1}}{(h-1)!}\right)^kdx.
\end{equation}

Now if we substitute (\ref{eq:ptnsgf}) into (\ref{eq:pgf}) we obtain the 
probability generating function of the generalized coupon collector's 
problem in the form
\begin{equation}
\label{eq:hcopies}
P_h(x)=\frac{1}{d^{h-2}}\int_0^{\infty}\left\{{xD-1\choose 
h-1}x^{h-1}e^{-td/x}\right\}\left(e^t-1-t-\dots-\frac{t^{h-1}}{(h-1)!}\right)^{d-1}dt.
\end{equation}
In the above, ${xD-1\choose h-1}$ is the differential operator that is 
defined by
\[{xD-1\choose 
h-1}f(x)=\frac{1}{(h-1)!}\left(x\frac{d}{dx}-1\right)\left(x\frac{d}{dx}-2\right)\dots 
\left(x\frac{d}{dx}-h\right)f(x).\]
However, it is easy to establish, by induction on $h$, the interesting fact 
that
\begin{equation}
\label{eq:weird}
{xD-1\choose h-1}x^{h-1}e^{-td/x}=\frac{(td)^{h-1}}{(h-1)!}e^{-td/x}.
\end{equation}
Hence we have proved the following evaluation.
\begin{theorem}
The probability generating function for the coupon collecting problem in 
which at least $h$ copies of each coupon are needed is given by
\begin{equation}
\label{eq:pgf2}
P_h(x)=\frac{d}{(h-1)!}\int_0^{\infty}t^{h-1}e^{-td/x}\left(e^t-1-t-\dots-\frac{t^{h-1}}{(h-1)!}\right)^{d-1}dt.
\end{equation}
\end{theorem}

\subsection{Two examples}
Let's look at the cases $h=1$, the classical case, and $h=2$, where we want 
to collect at least two specimens of each of the $d$ kinds of coupons.

If $h=1$ then (\ref{eq:pgf2}) takes the form
\[P_1(x)=d\int_0^{\infty}e^{-td/x}(e^t-1)^{d-1}dt.\]
If we expand the power of $(e^t-1)$ by the binomial theorem and integrate 
termwise we obtain
\[P_1(x)=xd\sum_{j=0}^{d-1}{d-1\choose j}\frac{(-1)^{d-1-j}}{d-jx},\]
which is precisely the partial fraction expansion of the classical 
generating function (\ref{eq:pgfh1}).

To see something new, let $h=2$. Then
\begin{equation}
P_2(x)=d\int_0^{\infty}te^{-td/x}\left(e^t-1-t\right)^{d-1}dt.
\end{equation}
Again, by termwise integration this can be made fairly explicit, but since 
the most interest attaches to the expectation, let's look at the average 
number of trials that are needed to collect at least two samples of each of 
$d$ coupons. This is $P_2'(1)$, which after some simplification takes the form
\begin{equation}
\label{eq:avgh2}
P_2'(1)=d^2\int_0^{\infty}\left(\frac{t^2}{e^t-1-t}\right)(1-(1+t)e^{-t})^ddt.
\end{equation}

 From this we can go in either of two directions, an exact evaluation or an 
asymptotic approximation. By termwise integration it is easy to obtain the 
following exact formula, which is a finite sum, for $\langle T\rangle_2$, 
the average number of trials needed to collect at least two of each of the 
$d$ kinds of coupons:
\begin{equation}
\label{eq:avgxact}
\langle T\rangle_2=d^2\sum_{m,j}(-1)^m{d-1\choose m}{m\choose 
j}\frac{(j+2)!}{(m+1)^{j+3}}
\end{equation}
For $d=2,3,4,5$ these are $2,11/2, 347/36,12259/864$. To facilitate 
comparison with the classical ($h=1$) case, we show below, for $1\le d\le 
10$, a table of the expected numbers of trials needed when $h=1,2$.

\[
\begin{array}{ccccccccccc}
\hfill d:&1&2&3&4&5&6&7&8&9&10\\
\langle 
T\rangle_1:&1.0000&3.0000&5.5000&8.3333&11.417&14.700&18.150&21.743&25.460&29.290\\
\langle 
T\rangle_2:&2.0000&5.5000&9.6389&14.189&19.041&24.134&29.425&34.885&40.492&46.230 

\end{array}
\]

\subsection{Asymptotics}
Now we investigate the asymptotic behavior of (\ref{eq:avgh2}), for large 
$d$, to compare it with the $d\log{d}$ behavior of the classical case where 
$h=1$.
\begin{theorem}
If there are $d$ different kinds of coupons, and if at each step we sample 
one of the $d$ kinds with uniform probability, let $\langle T\rangle_h$ 
denote the average number of samples that we must take until, for the first 
time, we have collected at least $h$ specimens of each of the $d$ kinds of 
coupons. Then for every $h\ge 1$, we have $\langle T\rangle_h\sim 
d\log{d}\quad (d\to\infty)$.
\end{theorem}
Consider first the case $h=2$. In (\ref{eq:avgh2}) we make the substitution
\begin{equation}
\label{eq:subst}
e^{-u}=1-(1+t)e^{-t},
\end{equation}
where $u$ is a new variable of integration. We then find that
\begin{equation}
\label{eq:newp}
P_2'(1)=d^2\int_0^{\infty}t(u)e^{-ud}du,
\end{equation}
where $t(u)$ is the inverse function of the substitution (\ref{eq:subst}), 
which is well defined since the right side of (\ref{eq:subst}) increases 
steadily from 0 to 1 as $t$ increases from 0 to $\infty$.

The main contribution to $P_2'(1)$ comes from values of $u$ near $u=0$, and 
when $u$ is near 0 we have
\[t(u)=-\log{u}+O(\log{\log{u}}).\]
Following the arguments in \cite{Wo}, sec. 2.2, we see that
$P_2'(1)$ of (\ref{eq:newp}) has the same asymptotic behavior as
\[d^2\int_0^c(-\log{u})e^{-ud}du,\qquad (0<c<1)\]
and in \cite{Wo} this is shown to be
\[\sim d^2\cdot \frac{\log{d}}{d}=d\log{d}.\]

Now we consider the asymptotic behavior of the expected number of trials 
for general values of $h$. {From} (\ref{eq:pgf2}) we see that this expected 
number of trials can be written in the form
\begin{equation}
\frac{d^2}{(h-1)!}\int_0^{\infty}\left\{\frac{t^h}{e^t-1-t-\dots 
-\frac{t^{h-1}}{(h-1)!}}\right\}\left\{1-\left(1+t+\frac{t^2}{2}+\dots 
+\frac{t^{h-1}}{(h-1)!}\right)e^{-t}\right\}^ddt
\end{equation}
Again we make the change of variable
\begin{equation}
\label{eq:newvar}
e^{-u}=1-\left(1+t+\frac{t^2}{2}+\dots +\frac{t^{h-1}}{(h-1)!}\right)e^{-t}
\end{equation}
in the integral, and it takes the remarkably simple form (compare 
(\ref{eq:newp}))
\[P_h'(1)=d^2\int_0^{\infty}t(u)e^{-ud}du,\]
where $t(u)$ is the inverse function of the substitution (\ref{eq:newvar}). 
Again the main contribution to the integral comes from small values of $u$, 
and when $u$ is small and positive we have
\[t(u)= -\log{u}+(h-1)\log{(-\log{u})}+\dots.\]
Using the method of sec. II.2 of \cite{Wo} once more, we find that
\begin{equation}
\label{eq:asym}
\langle T\rangle_h=d\log{d}+(h-1)d\log{\kern-2pt\log{d}}(1+o(1))\qquad 
(d\to\infty).
\end{equation}
We remark that in the case of $d=200$ coupons, the correct expected number 
of trials to obtain two of each coupon is 1614 trials, the approximation 
$d\log{d}$ is $1175$, and the approximation
  $d\log{d}+(h-1)d\log{\kern-2pt\log{d}}$ is $1393$, each rounded to the 
nearest integer.
\section{The number of singletons}
\label{sec:singletons}
In view of the asymptotics in the preceding section we realize that at the 
moment when a coupon collector sequence terminates with a complete 
collection, `most' coupons will have been collected more than once, and 
only `a few' will have been collected just once. We call a coupon that has 
been seen just once a \textit{singleton}. We will now look at the 
distribution of singletons.

In more detail, let $j$ be the number of singletons in a collecting 
sequence that terminates successfully at the $n$th step. We first want the 
joint distribution $f(n,j)$ of $n$ and $j$, i.e., the probability that a 
collecting sequence halts successfully at the $n$th step, and has exactly 
$j$ singletons at that moment. We claim that
\begin{equation}
\label{eq:joint}
f(n,j)=\frac{d!}{d^n}{n-1\choose j-1}\hstir{n-j}{d-j}{2}.
\end{equation}

Indeed, the last coupon to be collected can be chosen in $d$ ways, the 
other $j-1$ singleton coupons can be chosen in ${d-1\choose j-1}$ ways, and 
can be presented in an ordered sequence in $(j-1)!{d-1\choose j-1}$ ways. 
This ordered sequence can appear among the first $n-1$ trials in 
${n-1\choose j-1}$ ways, and the remaining $n-j$ trials constitute an 
ordered partition of $n-j$ elements into $d-j$ classes, no class having 
fewer than two elements, which can be chosen in $(d-j)!\hstir{n-j}{d-j}{2}$ 
ways. If we multiply these together and divide by $d^n$, the number of 
$n$-sequences, we obtain the result (\ref{eq:joint}) claimed above.

Next we compute the probability that a completed collecting sequence 
contains exactly $j$ singletons, whatever the length of the sequence may 
be. That is we find $F(j)=\sum_nf(n,j)$, where $f$ is given by 
(\ref{eq:joint}). We have, after using the generating function 
(\ref{eq:ptnsgf}),
\begin{eqnarray}
F(j)&=&\sum_n\frac{d!}{d^n}{n-1\choose j-1}\hstir{n-j}{d-j}{2}\\
&=&\frac{d!}{(d-j)!d^j}\left\{\int_0^{\infty}\left\{{t\frac{\partial}{\partial 
t}+j-1\choose 
j-1}\left(\frac{e^{-xt}}{t}\right)\right\}(e^x-1-x)^{d-j}dx\right\}_{t\rightarrow 
1/d}
\end{eqnarray}
But using the fact that, analogously to (\ref{eq:weird}), we have
\[{t\frac{\partial}{\partial t}+j-1\choose 
j-1}\left(\frac{e^{-xt}}{t}\right)=\frac{x^{j-1}}{(j-1)!t^j}e^{-x/t},\]
we can simplify the expression for $F(j)$ to
\begin{equation}
\label{eq:jprob}
F(j)=j{d\choose j}\int_0^{\infty}x^{j-1}(e^x-1-x)^{d-j}e^{-xd}dx,\qquad 
(j=1,2,3,\dots)
\end{equation}
which is the desired distribution of the number of singletons in a 
successfully terminated coupon collecting sequence.

Now if we multiply by $j$ and sum over $j$, we'll get the average number of 
singletons that appear in a completed collection of $d$ coupons. This is, 
after some termwise integration,
\[\bar{j}(d)=d\sum_m(-1)^m{d-2\choose m}\frac{d(m+1)+1}{(m+2)^2(m+1)}.\]
If we expand the summand in partial fractions, viz.,
\[\bar{j}(d)=d\sum_m(-1)^m{d-2\choose 
m}\left(\frac{1}{m+1}-\frac{1}{m+2}+\frac{d-1}{(m+2)^2}\right),\]
then each of the three sums indicated can be expressed in closed form, in 
two cases by using the identity
\begin{equation}
\label{eq:gould}
\sum_k(-1)^k{n\choose k}\frac{1}{x+k}=\frac{1}{x{x+n\choose n}},
\end{equation}
directly, with $x=1$ and $x=2$, and in the third case by differentiating 
(\ref{eq:gould}) w.r.t. $x$, and using the result with $x=2$. The identity 
(\ref{eq:gould}) is itself certified, after multiplying by the denominator 
on the right, by the WZ proof certificate $R(n,k)=k(x+k)/((n+1)(k-n-1))$.

What results is that $\bar{j}(d)=H_d$, the $d$th harmonic number. That is, 
\textit{the average number of singleton coupons in a completed collection 
sequence of $d$ coupons is the harmonic number $H_d$.}

\addcontentsline{toc}{section}{\hspace{.22in}References}

\end{document}